\documentclass{article}
\usepackage{amsfonts, amsbsy, amssymb}
\usepackage{graphicx}

\begin{document}


\newtheorem{theorem}{Theorem}
\newtheorem{proposition}{Proposition}[theorem]
\newtheorem{corollary}{Corollary}[theorem]
\newtheorem{lemma}{Lemma}[theorem]

\def \R{{\mathbb R}}
\def \T{{\mathbb T}}
\def \Z{{\mathbb Z}}
\def \eps{\varepsilon}
\title{Endpoint bounds for the non-isotropic Falconer distance
problem associated with lattice-like sets }
\author{A. Iosevich and M.
Rudnev}

\maketitle \begin{abstract} Let $S \subset {\mathbb R}^d$ be
contained in the unit ball. Let $\Delta(S)=\{||a-b||:a,b \in S\}$,
the Euclidean distance set of $S$. Falconer conjectured that the
$\Delta(S)$ has positive Lebesque measure if the Hausdorff dimension
of $S$ is greater than $\frac{d}{2}$. He also produced an example,
based on the integer lattice, showing that the exponent
$\frac{d}{2}$ cannot be improved. In this paper we prove the
Falconer distance conjecture for this class of sets based on the
integer lattice. In dimensions four and higher we attain the
endpoint by proving that the Lebesgue measure of the resulting
distance set is still positive if the Hausdorff dimension of $S$
equals $\frac{d}{2}$. In three dimensions we are off by a logarithm.

More generally, we consider $K$-distance sets
$\Delta_K(S)=\{{|a-b|}_K: a,b \in S\}$, where ${|\cdot|}_K$ is the
distance induced by a norm defined by a smooth symmetric convex body
$K$ whose boundary has everywhere non-vanishing Gaussian curvature.
We prove that our endpoint result still holds in this setting,
providing a further illustration of the role of curvature in this
class of problems.
\end{abstract}

\medskip
\noindent
{\bf Keywords:} Lattice point distribution, mean square estimates,\\
{\bf AMS subject classification} 42B

\subsection*{1. Introduction} Let $B$ denote the Euclidean unit ball.
The Falconer conjecture says that if the Hausdorff dimension of $S
\subset B\subset\R^d,\,d\geq2$ is greater than $\frac{d}{2}$, then
the Lebesgue measure of the distance set $\Delta(S)=\{\|a-b\|: a,b
\in S\}$ is positive. Here, and throughout the paper,
$\|x\|=\sqrt{x_1^2+\dots+x_d^2}$ is the Euclidean distance in
$\R^d$. See \cite{Fa}, \cite{Ma}, \cite{Bo}, \cite{Wo}, \cite{W},
\cite{Er}, and the references contained therein for the description
of this problem and progress over the years. The problem remains
open for every $d\geq2$. The best known results are due to Wolff in
${\mathbb R}^2$ and Erdogan in ${\mathbb R}^d$. They prove that the
Lebesgue measure of the distance set is positive if the Hausdorff
dimension of $S$ is greater than $\frac{d(d+2)}{2(d+1)}$.

Falconer showed that the exponent ${d\over2}$ cannot in general be
improved in the following sense. Let $\{q_i\}_{i\geq1}$ be a rapidly
growing sequence of positive integers, e.g with $q_1=2$ and
$q_{i+1}>q_i^i$. Let $S_i$ denote the union of Euclidean balls of
radius $q_i^{-\frac{d}{s}}$, for some $0<s<d$, centered at the
points of $\frac{1}{q_i} {\mathbb Z}^d \cap B$. Let $S=\cap_i S_i$.
Then (see e.g. \cite{F}) the Hausdorff dimension of $S$ is $s$. On
the other hand, the Lebesque measure
\begin{equation} |\Delta(S_i)| \leq {\rm const.}
\cdot q_i^{-\frac{d}{s}} \cdot q_i^2, \label{00}\end{equation} where
the first factor after the constant is the radius of each ball, and
the second factor is a trivial bound on the number of distances
determined by ${\mathbb Z}^d \cap {[0,q_i]}^d$. It follows that
$|\Delta(S)|=0$ if $s<\frac{d}{2}$. It is well known that in fact,
one can choose $\{q_i\}$ such that
\begin{equation} |\Delta(S_i)| = {\rm const.}\cdot q_i^{-\frac{d}{s}}
\cdot\left\{\begin{array}{ll} q_i^2, & d\geq3\\
\frac{q_i^2}{\sqrt{\log{q_i}}}, &d=2.\end{array}\right.
\label{01}\end{equation} Thus $|\Delta(S)|>0$  for $s\geq
\frac{d}{2}$ if $d\geq3$ and for $s> \frac{d}{2}$ if $d=2$.

As we mention in the abstract, we also study the case of more
general, "well-curved" distances. Let $K\subset \R^d$ be a strictly
convex symmetric body with a smooth boundary and the volume equal
to, say that of the Euclidean unit ball $B$. Let $|\cdot|_K$ be the
Minkowski functional of $K$, or the $K$-norm. Thus
\begin{equation}
\Delta_K(S)=\{|a-b|_K:\,a,b\in S\}.\label{02}
\end{equation}
Let $|\cdot |_{K^*}$ be the dual norm to $|\cdot |_{K}$, defined as
\begin{equation}
|x|_{K^*}=\sup_{y\in K}|x\cdot y|,\;\;\; K_*=\{x\in\R^d:
|x|_{K^*}\leq1\}.\label{a}
\end{equation}
The purpose of this paper is to prove that the Falconer conjecture
holds for the set $S$ constructed by Falconer, with respect to
$K$-distances, under the assumption that the boundary $\partial K$
of $K$ is $C^r$, for a large enough $r$ (we do not discuss what the
smallest possible $r$ could be) and has everywhere non-vanishing
Gaussian curvature. The fact that we are able to deal with any such
$K$ implies that as the basis for Falconer's construction one can
use any lattice, so considering $\Z^d$ as we do in the remainder of
this does not result in the loss of generality.

Furthermore, if one considers the Euclidean distance, we can extend
the scope of Falconer's construction to a homogeneous set
\begin{equation}
A=\left\{{|a|_K\over\|a\|}a:\,a\in\Z^d\right\},\;\mbox{ for some
}K,\label{ext}
\end{equation}
or equally well instead of the Euclidean norm $\|a\|$ in the
denominator we can have some other $|a|_{K'}$. Our main result is
the following.

\begin{theorem} Let $S$ be the Falconer set described above.
Suppose that $s\geq\frac{d}{2}$, and $d \ge 4$. Then the Lebesgue
measure of $\Delta_K(S)$ is positive, for any strictly convex $K$
with a smooth boundary. If $d=3$ the same conclusion holds with a
logarithmic loss in the sense to be made precise below.
\end{theorem}
The difficult part in Theorem 1 is the endpoint issue. Otherwise the
proof can be made somewhat shorter, using techniques developed by
M\"uller (\cite{Mu}) and Iosevich et al. (\cite{ISS}) for the
quantity
\begin{equation}{\cal E}(t)=
\# \{tK \cap {\mathbb Z}^d\}-t^d{\tt Vol}\,K.\label{aa}
\end{equation}

However, our main motivation in proving Theorem 1 is not the
quantity ${\cal E}$, but rather the Falconer distance conjecture,
and its discrete analog, the Erd\"os distance conjecture in the case
of homogeneous sets. An infinite discrete set $A\subset \R^d$ is
called homogeneous if all its elements are separated by some $c>0$,
while any cube of side length $C>c$ contains at least one element of
$A$. If $A_q=A\cap qK$ is a truncation of $A$, the conjecture is
that
\begin{equation}
\# \Delta_K(A_q) \geq C_\varepsilon q^{2-\varepsilon}. \label{aaa}
\end{equation}
This is a major open problem. See, for example, \cite{PS}, and the
references contained therein for the discussion of this problems and
the best known results. Theorem 2 below says that if $A$ is a
$d$-dimensional lattice, then
\begin{equation}
\# \Delta_K(A_q) \gtrsim \left\{
\begin{array}{ll}q^{2},&d\geq4,\\ q^{2}\log^{-2}q,
&d=3.\end{array}\right.\label{aaaa}
\end{equation} As we mention above,
this result is observed in \cite{ISSII} in the context of sharpness
examples for the lattice mean-square discrepancy estimates. The
symbol $\lesssim$ will further absorb constants depending only on
$K$ (and hence $d$). Also we write $a\gtrsim b$ if $b\lesssim a$ and
$a\approx b$ if both $a\lesssim b$ and $a\gtrsim b$. The symbol
$\sim$ will indicate proportionality, up to some constant $c(K)$.

Our approach is based on the formalism of distance measures, set up
by Mattila (\cite{Ma}). The distance measure $\nu(t)$, relative to
the set $A_q$ counts the number of its points in ${1\over q}$-thin
$K$-annuli of radius $t$, centered at points of $A_q$.
(Unfortunately, beyond Introduction, there will be several closely
related $\nu's$ which will have to carry extra identifications.) In
the case of a lattice it is enough to consider only the distances to
the origin, yet this is not essential. The distance measure
formalism has nothing to do with the lattice structure, and applies,
for example, to any homogeneous point set $A$. The $L^1$-norm
$\|\nu(t)\|_1$ simply counts the points in $qK$, and the main task
will be to estimate the $L^2$-norm $\|\nu(t)\|_2$. The latter can be
thought of as a special case of the estimation of the Mattila
integral. More precisely, Mattila (\cite{Ma}) proved a general
theorem for the Euclidean distance, which generalizes to
$K$-distances, (see \cite{AI}) as follows. Let $\mu$ be a finite
Borel measure supported on $S \subset B,$ Then if
\begin{equation} M(\mu)=\int_1^{\infty} {\left( \int_{\partial
K^{*}} {|\hat{\mu}(t \omega)|}^2 d\sigma_{K^*}(\omega) \right)}^2
t^{d-1}dt<\infty, \label{mi}\end{equation} where $d\sigma_{K^*}$ is
the Lebesgue measure on $\partial K^{*}$, the Lebesgue measure of
$\Delta_K(S)$ is positive. Falconer's example shows that there are
sets of dimension $s<{d\over2}$ such that the integral (\ref{mi})
diverges.

We shall write down an analog of the Mattila integral (\ref{mi}) for
a specific measure $\mu$ on the Falconer sets $S_i$, such that
$S=\cap S_{q_i}$ has dimension $s={d\over2}$. It will follow that
the integral (\ref{mi}), after appropriate scaling, can be bounded
by $O(1)$ in $d\geq4$ and by $\log^2 q_i$ in $d=3$.

\medskip
\noindent Unfortunately, our bounds on the Mattila integral can be
justified only if $A$ is a lattice. The reason is that in order to
get these bounds, we use a smooth approximation $E(t)$ of the
discrepancy ${\cal E}(t)$ as the auxiliary quantity, as $E(t)$
admits a well known representation via the Poisson summation
formula, and $E^2(t)$ looks rather similar to the integrand in the
Mattila integral.

Hence $L^2$ estimates for $\nu$ are closely related to the
$L^2$-estimates for $E$ in the lattice case. If one could establish
an analogue of this relation in the general homogeneous set context,
or find a proper substitute for $E$, one could then prove the
Erd\"os distance conjecture for homogeneous sets.

$L^2$-estimates for $E$ were obtained in \cite{Mu} and \cite{ISS}.
The technique in this paper is based on asymptotic methods, enabling
us to
\begin{itemize}
\item {\em essentially} dominate the desired $L^2$-estimate for
$E$ by a weighted $L^2$ estimate for $\nu$, see (\ref{nue}) below
\item which in turn can be dominated by another $L^2$-estimate for
$E$, see e.g. (\ref{vottak}) below.
\end{itemize}
The Poisson summation formula is absolutely essential for the
estimate in the first bullet above, but (seemingly) not for the
second one. As a result, one gets an $L^2$-estimate for the quantity
$E$, and a weighted $L^2$-estimate for the distance measure $\nu$.
It turns out that for $d\geq3$ the weighted $L^2$-estimate for $\nu$
implies a tight (modulo the logarithms in $d=3$) $L^2$-estimate for
$\nu$ itself, while for $d=2$ this is not the case. That is why our
method gives only a trivial estimate for $d=2$, where the case of a
general $K$ is an open problem.

We have italicized {\em essentially}, in the first bullet above,
because the asymptotic techniques used in  \cite{Mu} and \cite{ISS}
do not yield a clear cut relation (\ref{nue}) between the $L^2$
estimates for $\nu$ and $E$, due to the presence of various cut-off
functions, truncations, and related difficulties. This is precisely
the source of technical difficulties that one encounters in the
effort to attain the endpoint result claimed in Theorem 1. We
identify the estimate (\ref{nue}) as a key example of a technical
advantages of our approach. This approach yields the mean square
estimates for $E$ and ${\cal E}$ as a by-product. In addition,
throughout the proof a number of integral representations for the
distance measure $\nu$ and related quantities are obtained. These
identities extend to the case of general homogeneous sets and may
well prove to be of use in the further progress towards the Erd\"os
distance conjecture in this context.

\subsection*{2. Distance measure}
Let $\phi$ be a non-negative radial (radial henceforth means radial
with respect to the Euclidean metric) Schwartz class function, such
that $\int \phi(x)=1$, $\phi(x)=1$ inside the ball of some radius
and vanishes outside the ball of twice the radius. Let $q$ be a
large number, denote $\phi_q(x)=q^d\phi(qx)$ and $\Z^d_q=\Z^d\cap
qK\setminus\{0\}$.

For a function $f\in L^1(\R^d)\cap L^2(\R^d),$ let
\begin{equation}\hat f(\xi)=\int f(x)e^{-2\pi\iota \xi\cdot
x}dx\label{ft}\end{equation} define the usual Fourier transform. Let
\begin{equation} \mu_{q}(x)=\sum_{a\in \Z^d_q}\phi_q(x-a),\label{mu}
\end{equation}
be the smoothing of the counting measure on $\Z^d_q$, blurring each
point of $\Z^d_q$ into a speck of radius $\sim{1\over q}$ (actually
it is slightly less than ${1\over q}$ by the choice of $\phi$).
Clearly
\begin{equation}\hat\mu_{q}(\xi)= \sum_{a\in
\Z^d_q}\hat\phi(\xi/q)e^{-2\pi\iota
a\cdot\xi},\label{muhat}\end{equation} and the function $\hat\phi$
is radial.

For the construction of the sequence $S_i$ in the introduction, we
take a rapidly growing sequence $\{q_i\}$ and define $S_i$ to be the
scaling of $\Z_{q_i}^d$ into the unit cube. The resulting
intersection over $i\geq1$ has Hausdorff dimension ${d\over2}$ by
the argument (see e.g. \cite{F}) mentioned in the introduction.

For $t>0$ define
\begin{equation}\begin{array}{llllll}
\nu_{q,0}(t)&=&\int\omega_K(x/t)d\mu_{q}(x),\\ \hfill \\
N_{q,0}(t)&=&\int\Omega_K(x/t)d\mu_{q}(x)&=&\int_0^t d\nu_{q,0}.
\end{array} \label{nu}
\end{equation}
Above $\omega_K$ is the Lebesque measure on $\partial K$, $\Omega_K$
is the characteristic function of $K$. Note that in the first
integral $\mu_{q}$ is actually a Schwartz function, and $\omega_K$ -
a distribution.

Also define the volume discrepancy
\begin{equation} E_{q,0}(t)=N_{q,0}(t)-t^d{\tt Vol}K.
\label{e}\end{equation} Note that the quantity $E_{q,0}(t)$ thus
defined relates to ${\cal E}(t)$ in (\ref{aa}) for $t<q$ only. The
seemingly redundant ${\,}_0$ subscripts come from the fact that in
the sequel it turns out to be more convenient to work with the
weighted quantities
\begin{equation}
[\nu_q(t),\,N_q(t),
\,E_q(t)]\,=\,t^{1-d\over2}[\nu_{q,0},\,N_{q,0}(t), \,E_{q,0}(t)
].\label{wt}\end{equation} The quantity $\nu_{q,0}$ can be viewed as
the density of the measure $\mu_{q}$ on $K$-spheres of radius $t$,
centered at the origin. The primitive $N_{q,0}(t)$ counts the points
in $K$-balls of radius $t$. Removing the origin in $\Z^d_q$ is of no
consequence, yet it enables us to avoid consideration of some
trivial cases. Clearly
\begin{equation}\int_0^\infty\nu_{q,0}\;\sim \;q^d\label{l1}.
\end{equation}
By definition of the quantities $\mu_{q}$ and $\nu_{q,0}$, as we are
interested in the estimates only in terms of the order of magnitude
with respect to $q$, it is legitimate to sample the integrals
containing $\nu_{q,0}$ (as well as $E_{q,0}$ and other versions of
$\nu$ and $E$ to appear later) by Darboux sums with the step size
${1\over c_1 q}$, for some constant $c_1$.

Clearly $\nu_{q,0}$ vanishes for $t>q+q^{-1}$, while for
$t<q-q^{-1},$
\begin{equation}
{\nu_{q,0}\over q}\;\approx\;\Gamma(t,\delta), \label{gamma}
\end{equation} where $\Gamma(t,\delta)$ is the number of points of
$\Z^d$ in a $K$-annulus ${\cal A}(t,\delta)$ centered at the origin,
with radius $t$ and width $\delta$; in (\ref{gamma}) and further
$\delta$ will always be $\approx {1\over q}$. More precisely,
(\ref{gamma}) means that there exist uniform constants $c_2$ and
$c_3$, such that
\begin{equation}
\Gamma\left(t, {1\over c_2q}\right)\;\leq \;{\nu_{q,0}\over
q}\;\leq\;\Gamma\left(t, {c_3\over q}\right). \label{gamma1}
\end{equation}
Consider all consecutive annuli, indexed by $k\lesssim q^2$, of
width $\delta\sim {1\over q}$, whose $K$-radii $t_k$ go up to $q$.
Define the annulus standard deviation and $D_{\cal A}$ and the ball
mean square discrepancy $D_{K}$ as follows:
\begin{equation}\begin{array}{lllllll}
D_{{\cal A}_K}&=&\sqrt{{1\over q^2} \sum_k \Gamma^2(t_k,\delta) }
&\approx & \sqrt{{1\over q^3}\int_0^{q}
\nu^2_{q,0}(t)dt }, \\ \hfill \\
D_{K}&=& \sqrt{{1\over q^2}\sum_{k} E_{q,0}^2(t_k)} &\approx&
\sqrt{{1\over q}\int_{0}^{q} E^2_{q,0}(t)dt}.
\end{array}
\label{discr}
\end{equation}
We shall need the following estimate.
\begin{theorem} We have
\begin{equation}\begin{array}{llllll}
D_{{\cal A}_K},\,D_K& \lesssim  &q^{d-2},& &d\geq
4,\\ \hfill \\
D_{{\cal A}_K},\,D_K& \lesssim  &q\log q, && d=3.
\end{array}
\label{main}
\end{equation}
\end{theorem}
As far as the weighted quantity $\nu_q(t)$ is concerned, see
(\ref{wt}), the estimates of Theorem 2 boil down to establishing
\begin{equation}
\|\nu_{q}\|_2^2\;=\;\int_0^\infty \nu_{q}^2(t)dt \;\lesssim\;
\left\{\begin{array}{ll} q^d,&d\geq4,\\ q^3\log q,&d=3.
\end{array}\right.
\label{prove}
\end{equation}
Theorem 2 implies the following version of the Erd\"os distance
conjecture. See also \cite{ISSII} where this result is deduced in
the context of sharpness examples for

\begin{corollary} There are at least $\sim\,q^2$ distinct, ${1\over
q}$-separated $K$-distances from the origin to the points of
$\Z^d_q$, for $d\geq4$. If $d=3,$ $q^2$ changes to ${q^2\over\log^2
q}$. \label{cor}\end{corollary} We shall deduce Theorem 1 from
Corollary \ref{cor}.

\medskip\noindent
{\bf Proof of Corollary \ref{cor} and Theorem 1} Assume Theorem 2
for the moment. By Cauchy-Schwartz,
\begin{equation}q^{2d}\approx\left(\int_1^q
\nu_{q,0}dt\right)^2\leq |{\tt supp}\nu_{q,0}|\int_1^q
\nu_{q,0}^2(t)dt, \label{cs}\end{equation} where $|{\tt
supp}\nu_{q,0}|$ is the Lebesque measure of the support of
$\nu_{q,0}$. Now plug  the estimates (\ref{main}) in the right hand
side and get the lower bound $\gtrsim q$ for $d\geq4$ and $\gtrsim
{q\over\log^2 q}$ for $d=3$, for $|{\tt supp}\nu_{q,0}|$. The fact
that the distances are $\approx {1\over q}$-separated follows by
construction of $\nu_{q,0}$.

We now deduce Theorem 1 in the following way. Take a sequence
$\{q_i\}$ of values of $q$ as described in Introduction. Contract
$\Z^d_{q_i}$ into the unit ball to get the Falconer sets $S_i$
($S_i$ consists of points of $\Z^d_{q_i}$ scaled into the unit cube
and blurred into a smudge of radius of $\approx q^{-2}$). Then for
$d\geq 4$, the Lebesgue measure of each $\Delta_K(S_i)$ is bounded
uniformly away from zero, so the same conclusion must hold for the
distance set of $S=\bigcap S_i$. $\Box$

\medskip
\noindent Observe that in the same way as (\ref{cs}), yet for the
weighted quantity $\nu_{q}$ we always have
\begin{equation} \|\nu_{q}\|_2^2\;\gtrsim \;q^d.
\label{lwr}
\end{equation}

Let us now write up some representations for the quantities
$\nu_{q},N_{q}$, without using the Poisson formula, so they would
adapt to general distance measures on homogeneous sets.

First we apply Plahcherel to evaluate the integrals in (\ref{nu}).
Then for the weighted quantities (\ref{wt}) we get:
\begin{equation} \begin{array}{lll}\nu_{q}(t)&=&
t^{{d-1\over2}}\int\hat\phi(\xi/q)\hat\omega_K(t\xi)\sum_{a\in\Z^q_q}
e^{-2\pi\iota\,a\cdot\xi}d\xi,\\
\hfill \\
N_q(t)&=&
t^{{d+1\over2}}\int\hat\phi(\xi/q)\hat\Omega_K(t\xi)\sum_{a\in\Z^q_q}
e^{-2\pi\iota\,a\cdot\xi}d\xi.
\end{array}
\label{defs}\end{equation} Note that $\nu_{q}(t)$ extends as zero to
$t=0$, as well as the fact that $N_q(t)$ defined as it is not in
$L^2(\R_+)$ if $d=2$.

\subsubsection*{Euclidean case} We make a few remarks about the case
when $K$ is the Euclidean ball, with the notations $\omega_B$,
$\Omega_B$ for the surface and volume measure. In this case the
Fourier transform $\hat\omega_B(\xi)$ is radial, namely
$\hat\omega_B(\xi)\sim
\|\xi\|^{1-{d\over2}}J_{{d\over2}-1}(2\pi\|\xi\|)$, where further
$J_v$ denotes the Bessel function of order $v$. Let us skip the
factor of $2\pi$ in what follows. This can always be accomplished by
scaling. After writing the integral (\ref{defs}) for $\nu_{q}$ in
the spherical coordinates we have
\begin{equation}
\nu_{q}(t)\;\sim \;\sqrt{t}\int_0^\infty r J_{{d\over2}-1}(rt)
\psi(r/q)\sum_{a\in\Z^d_q} J_{{d\over2}-1}(r\|a\|)dr,
\label{ed}\end{equation} where henceforth
\begin{equation} \psi(r)\;= \;\hat\phi(\xi)|_{\|\xi\|=r},
\label{varphi}\end{equation} so $|\psi(r/q)|$ is asymptotically
smaller than any power of $r/q$.

After using Hankel's formula, see e.g. Watson's classic (\cite{Wa}),
\begin{equation}
\int_0^\infty t J_v(at)J_v(bt)dt\;=\;{\delta(a-b)\over a},
\label{hf}\end{equation}  for the Bessel functions of order $v\geq
0$, we would have then
\begin{equation} \begin{array}{lllllll}
\|\nu_{q}\|_2^2&=&\int_0^\infty \nu_{q}^2(t)dt &\sim & \int_0^\infty
r \psi^2(r/q) \sum_{a,b\in\Z^d_{q}}
{J_{{d\over2}-1}(r\|a\|)J_{{d\over2}-1}(r\|b\|)\over
(\|a\|\|b\|)^{{d\over2}-1}}dr
\\ \hfill \\
&&&\sim & \int_0^\infty r^{d-1} \psi^2(r/q) \sum_{a,b\in\Z^d_{q}}
\hat\omega_B(r\|a\|)\hat\omega_B(r\|b\|)dr.\end{array}
\label{mate}\end{equation} The representation (\ref{mate}) is a
particular case of the Mattila integral (\ref{mi}).

Note that due to Hankel's formula (\ref{hf}) with $v={d\over 2}-1$,
the expression (\ref{mate}) is in essence Parseval's identity for
the transformation
\begin{equation}
H[\nu_{q}](r)=\int_0^\infty\sqrt{r t}J_{{d\over2}-1}(r
t)\nu_{q}(t)dt\;\sim\; \sqrt{r} \psi(r/q)
\sum_{a\in\Z^d_{q}\setminus\{0\}} {J_{{d\over2}-1}(r\|a\|)\over
(\|a\|)^{{d\over2}-1}}.\label{ht}\end{equation} So, similarly to
(\ref{mate}) one can apply the Hankel formula with $v={d\over2}$ to
the quantity $N_q(t)$ and get for $d\geq3$
\begin{equation}
\begin{array}{lllllll}
q^{d+2}&\approx&\|N_q\|_2^2 &\sim & \int_0^\infty {1\over r}
\psi^2(r/q) \sum_{a,b\in\Z^d_{q}}
{J_{{d\over2}-1}(r\|a\|)J_{{d\over2}-1}(r\|b\|)\over
(\|a\|\|b\|)^{{d\over2}-1}}dr
\\ \hfill \\
&&&\sim & \int_0^\infty r^{d-3} \psi^2(r/q) \sum_{a,b\in\Z^d_{q}}
\hat\omega_B(r\|a\|)\hat\omega_B(r\|b\|)dr,\end{array}
\label{maten}\end{equation} It is easy to show that in order to get
the right order of magnitude $q^{d+2}$ for $\|N_q\|_2^2$ in the
latter integral, it suffices to restrict the domain of integration
to $(0,1)$, while the integral over $[1,\infty)$ should be $\approx
q^{-2}\|\nu_{q}^2\|_2^2$. For the lattice case studied, see the
ensuing Lemma \ref{lmmm}, the integral above, taken from $1$ to
infinity is closely related to the quantity $E_q(t)$.

However, for a general homogeneous set, even for the Euclidean
distance, there does not appear to be an analogous ``nice'' formula,
similar to (\ref{mate}), (\ref{maten}), for the quantity $E_q$.

\medskip \noindent
{\it Remark R1} In the next section we will show that in the case of
anisotropic distances $|\cdot|_K$, one has expressions similar to
(\ref{mate}) for $\|\nu_{q}\|_2^2$, using the asymptotics for the
Fourier transforms $\hat\omega_K$ and $\hat\Omega_K$. So far we have
not used the fact that we are dealing with a lattice. Hence, one can
squeeze $\Z^d_q$ in each direction $\xi\in S^{d-1}$ by the factor
$|\xi|_K$ and consider the Euclidean distances.

\subsubsection*{Anisotropic case}
To proceed, we need the following lemma on the asymptotics of the
Fourier transforms $\hat\omega_K$ and $\hat\Omega_K$. We do not
present a proof here, as for $\hat\Omega_K$ it can be found in
\cite{IR}, and the case of $\hat\omega_K$ can be treated in the same
way. One can also derive the expression for $\hat\omega_K$ from that
of $\hat\Omega_K$ using (\ref{defs}) and the fact that
$N_{q,0}(t)={d\nu_{q,0}\over dt}$. For more asymptotic expressions
of this kind see \cite{Hz}, \cite{So}.

\addtocounter{lemma}{1}\begin{lemma} For $\|\xi\|\leq 1$,
$\hat\omega_K(\xi),\, \hat\Omega_K(\xi)\approx 1$, otherwise
\begin{equation}\begin{array}{llllll}
\hat\omega_K(\xi)&=&\sum_{j=0}^1
u_j\left({\xi\over\|\xi\|}\right)J_{{d\over2}-1+j}
(c_4|\xi|_{K*})\|\xi\|^{1-{d\over2}-j}&+&O\left(\|\xi\|^{-{d+3\over2}}\right),
\\ \hfill \\
\hat\Omega_K(\xi)&=&\sum_{j=0}^1
U_j\left({\xi\over\|\xi\|}\right)J_{{d\over2}+j}
(c_4|\xi|_{K*})\|\xi\|^{-{d\over2}-j}&+&O\left(\|\xi\|^{-{d+5\over2}}
\right),
\end{array}
\label{ass}
\end{equation}
where the quantities $u_0, U_0$ are strictly positive and the
constant $c_4$ depends on $K$ only. \label{aspt}\end{lemma} Without
loss of generality, assume $c_4=1$ in the formulae (\ref{ass})
above. The sums in the asymptotic expansions have two terms, because
this is as many as we will have to analyze. Note that in the
Euclidean case, one has only the first term in the sum in the
expressions (\ref{ass}) and no remainder.

Lemma \ref{aspt} will be instrumental for our proofs. First let us
use it to elaborate on Mattila's integral.
\addtocounter{proposition}{2}\begin{proposition} For $d\geq2$,
\begin{equation} \begin{array}{lllllll}
\|\nu_{q}\|_2^2 &\approx& \int_0^\infty r^{d-1} \psi^2(r/q)
\sum_{a,b\in\Z_q^d} \hat \omega_{K^*}(ra) \hat\omega_{K^*}(rb) \,dr.
\end{array}\label{mat}
\end{equation}
\label{pp}
\end{proposition}
{\bf Proof} Let us see what the right-hand side of the first
expression in (\ref{mat}) is equal to by plugging in the asymptotic
expansion for $\hat \omega_{K^*}$ from (\ref{ass}). Keep in mind
that $K^{**}=K$. Given $(a,b),$ it suffices to consider three cases,
as far as the three-term expansions in (\ref{ass}) are concerned:
the leading terms for both $a$ and $b$, the leading term for $a$ and
the second term for $b$, and finally the leading term for $a$ and
the remainder for $b$. As (\ref{ass}) has three terms, there will be
three key estimates in the proof, and the same will apply to proofs
of Theorem 3 and Lemma \ref{lmmm} in the sequel.

{\bf 1.} As far as the leading term in the sum (\ref{ass}) is
concerned, the fact that its contribution in the integral
(\ref{mat}) is $\approx \|\nu_{q}\|_2^2$ follows immediately from
Remark R1. Let us also verify it directly however, because the same
argument will be used further for Theorem \ref{dtm} and Lemma
\ref{lmmm}. Given $(a,b)$, take the product of the leading terms in
the sum (\ref{ass}). Then, see (\ref{mate}), what we get is
proportional to the following integral over $\R^d$:
\begin{equation}\begin{array}{rr}
u_0(a/\|a\|)u_0(b/\|b\|)\int \hat\phi^2(\xi/q) \hat
\omega_{B}(|a|_K\xi) \hat\omega_{B}(|b|_K\xi)d\xi \\ \hfill
\\ \approx \; \int [\phi_q*(\omega_B\circ {|a|^{-1}_K})](x)\cdot
[\phi_q*(\omega_B\circ {|b|^{-1}_K})](x)dx,\end{array}
\label{an}\end{equation} where $\omega_B\circ
{|a|_K^{-1}}(x)=|a|^{1-d}_K \omega_B(|a|_K^{-1}x).$ The integrand in
the right hand side of (\ref{an}) is roughly the product of
characteristic functions of concentric Euclidean annuli, of radii
$|a|_K$ and $|b|_K$, in particular it vanishes if
$||a|_K-|b|_K|>{1\over q}$, and is proportional to $|a|_K^{d-1}$ in
case $|a|_K=|b|_K$. Thus the average value of $\omega_B\circ
{|a|^{-1}_K}$ across the Euclidean annulus of radius $|a|_K$ and
width of $\approx {1\over q}$ is proportional to $q$. So we get
\begin{equation}
\sum_{a,b\in\Z^d_q} u_0(a/\|a\|) u_0(b/\|b\|) \int_0^\infty
r\psi^2(r/q) {J_{{d\over2}-1}(|a|_Kr)J_{{d\over2}-1}(|b|_Kr)\over
(|a|_K|b|_K)^{{d\over2}-1}}dr \;\approx \; q\sum_{k=1}^{\sim q^2}
{\Gamma^2(r_k,\delta)\over r_k^{d-1}}\;\approx\;\int_{0}^\infty
\nu_{q}^2(t)dt, \label{prc}\end{equation} where the sum above is
taken over consecutive $K$-annuli of radius $r_k$ and fixed width
$\delta$ of $\approx {1\over q}$. Recall that $\Gamma(r_k,\delta)$
is the cardinality of the intersection of $\Z^d$ with such an
annulus, see (\ref{gamma}).

{\bf 2.} Now let us take the first term in the sum (\ref{ass}) for
$a$ and the second one for $b$, plugging them in the right-hand side
of (\ref{mat}). Note that merely using the leading order asymptotics
in this case results in a superfluous factor $\int_1^\infty r^{-1}
\psi^2(r/q)dr \approx \log q$.

So, given $(a,b)$ we have, similarly to (\ref{an}), and omitting
uniform constants:
\begin{equation} \begin{array}{lll}
\int_0^\infty \psi^2(r/q)
{J_{{d\over2}-1}(|a|_Kr)J_{d\over2}(|b|_Kr)\over
|a|^{{d\over2}-1}_K|b|_K^{d\over2}}dr &\sim & \int \hat\phi^2(\xi/q)
\hat \omega_{B}(|a|_K\xi)
\hat\Omega_{B}(|b|_K\xi)d\xi  \\ \hfill \\
& = & \int [\phi_q*(\omega_B\circ {|a|^{-1}_K})](x) \cdot
[\phi_q*(\Omega_B\circ
{|b|^{-1}_K})](x)dx \\ \hfill \\
& \approx &\left\{\begin{array}{ll} {1\over |a|^{d-1}_K |b|^{d}_K}
|a|^{d-1}_K,& |a|_K\leq |b|_K+{1\over q}
,\\0&\mbox{otherwise.}\end{array}\right.
\end{array}\label{sec}\end{equation}
Hence, summing over $a,b\in\Z^d_q,$ we get \begin{equation} \approx
\;\sum_{b\in\Z_q^d\setminus\{0\}} |b|_K^{-d} \sum_{|a|_K\leq
|b|_K}1\;\approx\; q^d.\label{seca}
\end{equation}

{\bf 3.} Finally, as far as the remainder in (\ref{ass}) is
concerned, let us rewrite the integral in (\ref{mat}) as
$\sum_{a,b\in\Z^d_q} I_{a,b}$ and notice that without loss of
generality one can assume $|a|_K\geq|b|_K$. Then partition
\begin{equation} \sum_{a,b\in\Z^d_q,\,|b|_K\leq|a|_K}
I_{a,b}\;=\;\sum_{a,b\in\Z^d_q,\,|b|_K\leq|a|_K}\left(\int_0^{|a|^{-1}_K}
+ \int_{|a|^{-1}_K}^{|b|^{-1}_K}+\int_{|b|^{-1}_K}^\infty\right).
\label{split}\end{equation} The first piece - plug $1$ for the
$a$-term and $1$ for the $b$-term - can be bounded via
\begin{equation}
\sum_{a\in\Z^q_d}\; \sum_{b\in\Z^q_d,\,|b|_K\leq |a|_K}
\int_0^{|a|_K^{-1}}r^{d-1}dr\;\approx \;\sum_{a\in\Z^q_d} |a|_K^{-d}
\sum_{b\in\Z^q_d,\,|b|_K\leq |a|_K}1\;\approx\;
\sum_{a\in\Z^q_d}1\;\approx \;q^d. \label{odin} \end{equation} For
the second piece - plug $1$ for the $b$-term and zero order
asymptotics $(|a|_K r)^{-{d-1\over2}}$ for the $a$-term - we have a
bound
\begin{equation}
\sum_{a\in\Z^q_d} \;\sum_{b\in\Z^q_d,\,|b|_K\leq |a|_K}
\int_{|a|^{-1}_K}^{|b|^{-1}_K} |a|^{1-d\over2}_K
r^{d-1\over2}dt\;\approx\; \sum_{a\in\Z^q_d} |a|^{1-d\over2}_K
\sum_{b\in\Z^q_d,\,|b|_K\leq |a|_K} |b|_K^{-{d+1\over2}}\;\approx\;
\sum_{a\in\Z^q_d}1\;\approx\; q^d. \label{dva} \end{equation}
Finally, for the third piece plug $(|a|_K r)^{-{d-1\over 2}}$ for
the $a$-term terms and $(|b|_K r)^{-{d+3\over 2}}$, for the $b$-term
to get
\begin{equation}
\left(\sum_{a,b\in\Z_q^d\setminus\{0\}} |a|_K^{-{d-1\over2}}
|b|_K^{-{d+3\over2}} \right)\cdot\int_1^\infty
\psi^2(r/q)t^{-2}dt\;\lesssim \;q^{d-1}. \label{tri}
\end{equation} Hence, the upper estimates in the second and the
third cases match the lower bound (\ref{lwr}). This completes the
proof of Proposition \ref{pp} $\Box$.

\subsection*{3. Poisson formula}
As the set $\Z_q^d$ has beet taken off the integer lattice, the
values of $\nu_{q}(t),N_q(t),E_q(t)$ can be computed directly,
rather than via (\ref{defs}), using the Poisson summation formula.
It gives "nice" expressions for these quantities on the $t$-side,
rather than on the Hankel transform side, see (\ref{ht}). For
example, for ${1\over q}<t<q-{1\over q}$, one has
\begin{equation}
\int\omega_K(x/t)\sum_{a\in\Z^d_q}\phi_q(x-a)dx=\int\omega_K(x/t)
\sum_{a\in\Z^d}\phi_q(x-a)dx=\sum_{b\in
t^{-1}\Z^d}[\omega*\phi_{qt}](b),\label{conv}\end{equation} where
$\phi_{qt}(x)=\phi_q(tx).$ Applying the Poisson summation formula to
the convolution in the square brackets, not forgetting the scaling
(\ref{wt}) and doing the same thing for the quantities $N_q,E_q$
yields, for ${1\over q}<t<q-{1\over q}$
\begin{equation} \begin{array}{lllllllll}\nu_{q}(t)&\sim&
t^{{d-1\over2}}\sum_{a\in\Z^q}\hat\phi(a/q)\hat\omega_K(ta)&\equiv&\nu(t),
\\ \hfill \\
N_q(t)&\sim&
t^{{d+1\over2}}\sum_{a\in\Z^q}\hat\phi(a/q)\hat\Omega_K(ta)&\equiv&N(t),
\\ \hfill \\
E_q(t)&\sim&
t^{{d+1\over2}}\sum_{a\in\Z^q\setminus\{0\}}\hat\phi(a/q)\hat\Omega_K(ta)&\equiv&E(t).
\end{array}
\label{poisson}
\end{equation}
Note however that the quantities in the right hand side of
(\ref{poisson}) are unbounded as $t\rightarrow\infty$, besides the
summation is carried over the whole integer lattice. Still, there is
a considerable resemblance between the expression for $\nu(t)$ and
the square root of the integrand in (\ref{mat}). Let us introduce
the $L^2(\R_+)$ quantities
\begin{equation}
[\tilde\nu(t),\,\tilde{N}(t),\,\tilde{E}(t)]\;=\;\psi(t/q)[\nu(t),\,N(t),\,E(t)].
\label{tld}
\end{equation}
Clearly
\begin{equation}
\|\nu_q\|_2^2\;\lesssim\;\|\tilde\nu\|_2^2,\;\;\;\int_0^q
E_q^2(t)dt\;\lesssim \|\tilde{E}\|_2^2. \label{nuks}\end{equation}

\noindent {\it Remark R2} Note that in the definition (\ref{tld}) of
$\tilde\nu$ there is no harm restricting the summation to
$\Z^d\setminus\{0\},$ which will be done further. Indeed, $a=0$
results in a regular term $t^{d-1\over 2}$, and it is easy to check
that the contribution of this term into $\|\tilde\nu\|_2^2$ is
$O(q^d)$, cf. (\ref{lwr}).

\medskip
\noindent Now define the quantities $\nu_*(t),\,\tilde\nu_*(t)$ in
the same way as $\nu(t), \tilde\nu(t)$, has been defined, cf.
(\ref{nu}), (\ref{wt}), (\ref{tld}) but only with respect to the
dual body $K^*$ rather than $K$, do the same for $E$. Also with
respect to $K^*$, define the notation $\Gamma_*$, cf. (\ref{gamma}).
The next theorem is a bit of a red herring, as it is not used to
prove Theorem 1, yet is interesting in its own right. Of course,
this theorem is obvious in the Euclidean or ellipsoidal cases!

\begin{theorem}
For $d\geq2$,
\begin{equation}
\|\tilde\nu\|_2 \;\approx \; \|\tilde\nu_*\|_2. \label{dty}
\end{equation}
\label{dtm}\end{theorem} {\bf Proof} The proof follows the same
three steps as does the proof of Proposition \ref{pp}, adding to it
the decay of $\psi$ only. Let us use the definitions
(\ref{poisson}), (\ref{tld}) and the asymptotics (\ref{ass}) to
estimate the left-hand side.

{\bf 1.} For the principal term, skipping $u_0>0$, we get an
analogue of (\ref{prc}):
\begin{equation}\begin{array}{llllll}
\sum_{a,b\in\Z^d\setminus\{0\}} \hat\phi(a/q)\hat\phi(b/q)
\int_0^\infty r\psi^2(r/q)
{J_{{d\over2}-1}(|a|_Kr)J_{{d\over2}-1}(|b|_Kr)\over
(|a|_K|b|_K)^{{d\over2}-1}}dr &\approx &
q\sum_{k=1}^\infty{\Gamma_*^2(r_k,\delta)\over
r_k^{d-1}}\psi^2(r_k/q)\\
\hfill \\
&\approx&  \|\tilde\nu_*\|_2^2.\end{array}
\label{prcc}\end{equation}

{\bf 2,3.} For the rest of the estimates, one has chase through
(\ref{sec}--\ref{tri}), for the only difference that the summations,
say in $a$ therein will have been extended to
$\|a\|\rightarrow\infty$, and weighted by $(1+\|a\|/q)^{-N}$, for
some large $N$, reflecting the decay of $\psi$. Such an extension
can be easily verified to be of no consequence for the order of
these estimates. For instance (\ref{seca}) turns into a bound
\begin{equation}
\int_1^\infty {r^{d-1}\over (1+r/q)^N}dr\;=\;O(q^d), \label{chg}
\end{equation}
the same non-consequential changes happen to
(\ref{odin})--(\ref{tri}). $\Box$

The next lemma makes precise the statement in the first bullet in
the Introduction.
\begin{lemma} We have for $d\geq2$:
\begin{equation}
\|\tilde{E}\|^2_2\;\lesssim\; \int_0^\infty {\tilde\nu_*^2(t)\over
1+t^2} dt \;+ \;R(q),\label{nue}\end{equation} where
$R(q)=O(q^{d-2})$ in $d\geq3$ and $O(\log q)$ in $d=2$.
\label{lmmm}\end{lemma} {\bf Proof} The proof follows the same
pattern as proofs of Proposition \ref{pp} and Theorem \ref{dtm}. We
use the definitions (\ref{poisson}) and (\ref{tld}) and plug in the
asymptotics (\ref{ass}).

{\bf 1.} For the principal terms, skipping $U_0>0$, we get
\begin{equation}
\sum_{a,b\in\Z^d\setminus\{0\}}
\hat\phi(a/q)\hat\phi(b/q)\int_0^\infty t \psi^2(t/q)
{J_{d\over2}(|a|_{K*}t)J_{d\over2}(|b|_{K*}t)\over
(|a|_{K*}b_{K*})^{d\over2}}\,dt. \label{dura}
\end{equation}
The integral in (\ref{dura}), given $(a,b)$ can be rewritten as an
integral over $\R^d$:
\begin{equation}\begin{array}{r}
\int [(\Omega_B\circ|a|^{-1}_{K*})  * \phi_q ]^{\widehat{\,}}(\xi)
\xi \cdot [(\Omega_B\circ|b|^{-1}_{K*})
* \phi_q ]^{\widehat{\,}}(\xi)\xi d\xi \\ \hfill \\ =\;
\int \nabla_x[(\Omega_B\circ|a|^{-1}_{K*})  * \phi_q ](x)
\cdot\nabla_x[(\Omega_B\circ|b|^{-1}_{K*})
* \phi_q ](x)dx.\end{array}
\label{duraa} \end{equation} where, cf. (\ref{an}),
$\Omega_B\circ|a|^{-1}_{K*}(x)=|a|^{-d}_{K*}\Omega_B(|a|^{-1}_{K*}x)$.
The integral in the right-hand side of (\ref{duraa}) is clearly zero
if $||a|_{K^*}-|b|_{K^*}|>{1\over q}$, while if
$|a|_{K^*}=|b|_{K^*}$, it is $O(|a|^{d-1}_{K^*})$. Hence,
(\ref{dura}) is
\begin{equation}
\approx \;q\sum_{k=1}^\infty{\Gamma_*^2(r_k,\delta)\over r_k^2
r_k^{d-1}}\psi^2(r_k/q)\\ \hfill \\
\;\approx\; \int_0^\infty {\tilde\nu_*^2(t)\over 1+t^2} dt,
\label{escheraz}\end{equation} see also (\ref{prc}) and
(\ref{prcc}).

{\bf 2.} For the principal $a$-term and second $b$-term in the
asymptotics (\ref{ass}), throwing away uniform constants, we get

\begin{equation}
\sum_{a,b\in\Z^d\setminus\{0\}}
\hat\phi(a/q)\hat\phi(b/q)\int_0^\infty \psi^2(t/q)
{J_{d\over2}(|a|_{K*}t)J_{{d\over2}+1}(|b|_{K*}t)\over
|a|^{d\over2}_{K*} b_{K*}^{{d\over2}+1}}\,dt. \label{durab}
\end{equation}
Clearly (\ref{durab}) is reminiscent of (\ref{sec}), only in
dimension $d+1$. Let $w_B,W_B$ be the Lebesgue measure on $S^d$ and
the characteristic function of the Euclidean unit ball in
$\R^{d+1}$, respectively. Let $(y,\zeta)\in\R^{d+1}\times\R^{d+1},$
let the radial cutoff function $\varphi$ be defined in the same way
as $\phi$, only in dimension $d+1$.

Denote $\varphi_q(y)=q^{d+1}\varphi(qy)$,
$w_B\circ|a|_{K*}^{-1}(y)=|a|_{K*}^{-d}w_B(|a|_{K*}^{-1}y),\;
W_B\circ|a|_{K*}^{-1}(y)=|a|_{K*}^{-d-1}W_B(|a|_{K*}^{-1}y)$. Then
the integral right hand side of (\ref{durab}) can be rewritten as
follows:
\begin{equation}
\begin{array}{lll} \int_0^\infty \psi^2(t/q)
{J_{d\over2}(|a|_{K*}t)J_{{d\over2}+1}(|b|_{K*}t)\over
|a|^{d\over2}_{K*} b_{K*}^{{d\over2}+1}}\,dt&\sim& \int
[(w_B\circ|a|_{K*}^{-1})*\varphi_q]^{\widehat{\,}}(\zeta)\cdot
\left([(W_B\circ|a|_{K*}^{-1})*\varphi_q
]^{\widehat{\,}}(\zeta)\,\|\zeta\|\right)\,d\zeta \\ \hfill \\
&\approx & \int [(w_B\circ|a|_{K*}^{-1})*\varphi_q](y) \cdot
\|\nabla_y [(W_B\circ|a|_{K*}^{-1})*\varphi_q](y)\|\,dy.
\end{array}
\label{durabl}
\end{equation}
Thus the integral vanishes if $||a|_{K^*}-|b|_{K^*}|>{1\over q}$ and
is approximately ${1\over |b|_{K^*}^{d+1}}$ if
$|a|_{K^*}=|b|_{K^*}$. Summing over $(a,b)$ and taking absolute
values, we get the bound (\ref{escheraz}) once more.

{\bf 3.} We deal with the remainder in (\ref{ass}) in the same way
as it was done in Proposition \ref{pp} and Theorem \ref{dtm}, to
show that its contribution is $O(q^d)$. The demonstration is
routine. One simply chases through (\ref{split}--\ref{tri}),
modifying the exponents involved in the obvious way and extending
summations to infinity, yet weighting them by $(1+|a|_{K^*}/q)^{-N},
\,(1+|b|_{K^*}/q)^{-N},$ which is not consequential for the order of
magnitude of these estimates in $q$. Note that an extra $\log q$
gets picked up in $d=2$ in (\ref{odin}) where $|a|_K^{-d}$ changes
to $|a|_K^{-d-2}$. $\Box$

\medskip
\noindent {\it Remark R3} It is clear that estimating the
right-hand-side of (\ref{nue}) which is essentially an
$L^2$-estimate for ${\nu\over t}$ does not suffice to get a sharp
estimate for $\nu$, when it grows in average slower than $\sqrt{t}$
as $t\rightarrow\infty$. That is why we cannot prove Theorem 1 for
$d=2$. The necessity of the logarithmic factor in $d=3$ in the
estimate for $\nu$ is also questionable.

\subsection*{4. Proof of Theorem 2}
Theorem 2 will follow immediately from the bound (\ref{nue}) of
Lemma \ref{lmmm} and the following lemma, which somewhat generalizes
the results of \cite{ISS}. \addtocounter{theorem}{1}
\begin{lemma}\label{ebound} We have the following bound:
\begin{equation}\label{ebalda}
\|\tilde{E}\|^2_2\;\lesssim \;b_d(q),\;\mbox{ where
} \;b_d(q)\;=\;\left\{\begin{array}{lll}q^{d-2},& d\geq4,\\
q\log^2 q, &d=3, \\ q, & d=2.\end{array}\right.\end{equation}
\end{lemma}
{\bf Proof} There is no harm changing in(\ref{nue}) the lower limit
of integration to $1$ and $1+t^2$ in the denominator to $t^2$.

By definition of $\tilde\nu$, for any $t_l, t_u$, with
$t_u-t_l\gg{1\over q}$ and a small enough $\delta\sim{1\over q}$, we
have the representation of the integral as a Darboux sum:
\begin{equation}\int_{t_l}^{t_u}{\tilde\nu^2_*(t)\over t^2}dt\;
\approx \;\sum_{k} {\nu^2_*(t)|_{t\in I_k}\over t_k^2}\psi(t_k/q)\,
\delta,\label{dsum}
\end{equation}
where the intervals $I_k=[t_k,t_{k+1})$ of length $\delta$ partition
$[t_l,t_u)$ and choice of $t\in[t_k,t_{k+1})$ is arbitrary.

Then one can always choose $t$ inside each interval $I_k$ in such a
way that \begin{equation}\label{sbl}\nu_*(t)\;\lesssim\;
\max\left[t^{d-1\over2},\,q|E_*|(t)\right].\end{equation} This claim
quantizes the second bullet assertion in Introduction.

Indeed, the first term inside the above maximum corresponds to the
case of the existence of $t\in I_k$ such that $\nu_*(t)\lesssim
t^{d-1\over2}$. Otherwise, let us use the fact that ${dE_*(t)\over
dt} \approx \nu_*(t) + O\left(t^{d-1\over2}\right)$, and if
$\nu_*(t)\gtrsim t^{d-1\over2}$, the $O\left(t^{d-1\over2}\right)$
term can be omitted. Then $|E_*(t)|\gtrsim \int_{t_0}^t
\nu_*(\tau)d\tau$, where at $t_0, \,|E_*|$ has its absolute minimum
in $I_k$. Which implies that $q\sup_{I_k}| E_*(t)|\gtrsim \inf_{I_k}
\nu_*(t)$ in this case.

Note that due to (\ref{lwr}) all the ``regular'' terms
$O\left(t^{d-1\over2}\right)$ that pop up further would a-priori
result in (\ref{ebalda}), and in fact  stronger inequalities for
$d=2,3$.

Furthermore by (\ref{sbl}) \begin{equation} \int_1^\infty
{\tilde\nu^2_*(t)\over t^2}dt \;\lesssim \;\int_1^\infty
t^{d-3}\psi(t/q) dt\;+\; \int_{\cal I}{\tilde\nu^2_*(t)\over t^2}dt,
\label{aha}
\end{equation} where
\begin{equation}
{\cal I}\;=\;\{t: \nu_*(t)\leq c_5 q|E_*(t)|\}, \label{ahaha}
\end{equation}
for some $c_5$. The first integral in (\ref{aha}) bounded via
$q^{d-2}$ for $d\geq 3$ and $\log q$ for $d=2$.

Let us turn to the second integral in (\ref{aha}). Evaluation of
this integral goes along the lines of \cite{Mu}, \cite{ISS}.
Clearly, in order to get the upper bound, the integral can be
extended from ${\cal I}$ to $\R_+$, under the assumption that
$\nu_*(t)\leq c_5 q|E_*(t)|$ everywhere (note that ${\cal I}$ can be
represented as the union of intervals of length not smaller than
$\approx {1\over q}$ each). Under this assumption, we write out a
dyadic decomposition:
\begin{equation}
\int_1^\infty {\tilde\nu^2_*(t)\over 1+t^2}dt\;\approx\;
\sum_{k=0}^{\infty} 2^{-2k}\int_{2^k}^{2^{k+1}}\tilde\nu^2_*(t)dt
\;\lesssim\; q\sum_{k=0}^{\infty} |\psi(2^k/q)| 2^{{d+1\over4}k-2k}
\sqrt{ \int_{2^k}^{2^{k+1}}|\tilde{E}_*|^2\nu_*(t)dt
},\label{ddc}\end{equation} To get the right-hand side we have
applied Cauchy-Schwartz and used the fact that in an annulus of
width $2^k$ the integral of $\nu_*$ is $O(2^{k{d+1\over2}})$ -
recall the scaling (\ref{wt}).

Furthermore, using the fact that ${dE(t)\over dt} \approx \nu_*(t) +
O(t^{d-1\over2})$, we have
\begin{equation}
\int_{2^k}^{2^{k+1}}|\tilde{E}_*|^2\nu_*(t)dt \;\lesssim\;
|\psi(2^k/q)|\left(|E_*(2^k)|^3+ |E_*(2^{k+1})|^3+
2^{k{d-1\over2}}\int_{{2^k}}^{2^{k+1}}{E}_*^2(t)dt\right).
\end{equation}
The cubic terms in brackets are bounded as
\begin{equation}\label{ept}O\left[2^{3k\left({d-3\over2}+
{2\over d+1}\right)}+2^{3k{d-1\over 2}}q^{-3}\right],
\end{equation} which follows from the well known, see e.g.
Landau's classic (\cite{La}), $L^\infty$ estimate
\begin{equation}\label{linf} |E_0(t)|\;\lesssim\; t^{d-2+{2\over
d-1}}+ q^{-1}t^{d-1},\end{equation}, where
$E_0(t)=t^{{d-1\over2}}E(t)$, in view of the scaling (\ref{wt}). It
is a routine calculation to show using the decay of $\psi$ that the
contribution of these terms into (\ref{ddc}) is well in compliance
with (\ref{ebalda}).

Hence we are left with
\begin{equation}
\int_0^\infty \tilde{E}^2(t)dt \;\lesssim \;b_d(q)\, + \,q
\sum_{k=0}^{\infty} 2^{k\left({d\over2}-2\right)}|\psi(2^k/q)|
\sqrt{\int_{2^k}^{2^{k+1}} \tilde{E}_*^2(t)dt}. \label{vottak}
\end{equation}
Assuming that the sum above is $\gtrsim b_d(q)$, consider the case
$d\geq4$ first. Then, as clearly
\begin{equation}\sum_{k=0}^{\infty}
2^{k\left({d\over2}-2\right)}|\psi(2^k/q)|\;\lesssim\;
q^{{d\over2}-2},\;\mbox{ for }\;d\geq4, \label{hm}\end{equation} we
have
\begin{equation}
\int_0^\infty \tilde{E}^2(t)dt \;\lesssim \;q^{{d\over2}-1}
\sqrt{\int_{0}^{\infty} \tilde{E}_*^2(t)dt}, \label{vottaktak}
\end{equation}
and it follows that
$\|\tilde{E}\|_2^2,\,\|\tilde{E}_*\|_2^2\;\lesssim
\;b_d(q)\,=\,q^{d-2},\;d\geq4$, as one can certainly swap the
subscript ${\,}_*$ to the left-hand side.

The case $d=2,3$ requires some extra consideration, see \cite{ISS},
which we have adopted from the latter reference for the sake of
completeness. Recall that the quantity $E$ has been defined with
respect to the parameter $q$, where ${1\over q}$ is the
characteristic size of a speck, each point has been blown up to. To
reflect this fact, let us further write
$E=E^{(q)},\,\tilde{E}=\tilde{E}^{(q)}.$ It is easy to verify by
definition of $E$ that for $t\lesssim\bar{q}\lesssim q$, one has
\begin{equation}\label{ustal}|E^{(q)}|(t)\;\lesssim\;
|E^{(\bar{q})}|(t)+O(t^{d-1\over2}\bar{q}^{-1}).\end{equation}
Qualitatively this means that as the radial density of lattice
points increases with the radius, so the points near the origin can
be blown up into larger specks, without a risk of miscounting.

We rewrite (\ref{vottak}) as follows:
\begin{equation}
\int_0^\infty  |\tilde{E}^{(q)}|^2(t)dt \;\lesssim \;b_d(q)\, +\,
q\sup_k\left(\sqrt{{\int_{0}^{2^{k+1}}
|\tilde{E}_*^{(q)}|^2(t)dt\over b_d(2^{k+1})}}\right)
\sum_{k=0}^{\infty}
2^{k\left({d\over2}-2\right)}|\psi(2^k/q)|\sqrt{b_d(2^{k+1})},
\label{pochti}
\end{equation}
therefore, evaluating the sum we get
\begin{equation}
{\int_0^\infty  |\tilde{E}^{(q)}|^2 dt \over b_d(q)} \;\lesssim
\;1\,+\, \sup_k\left(\sqrt{{\int_{0}^{2^{k+1}}
|\tilde{E}_*^{(q)}|^2dt\over b_d(2^{k+1})}}\right).
\label{sovsempochti}
\end{equation}
The supremum above is at most some power of $q$. It should also be
achieved for some finite $k$, because of the decay, built into the
quantity $\tilde{E}$, due to the presence of the cutoff $\psi$.
Consider then such $\bar{k}$, when
\begin{equation}\label{nuda}
m(q)\;=\;\max\,\sup\left(
{\int_{0}^{2^{k+1}}|\tilde{E}^{(q)}|^2dt\over
b_d(2^{k+1})},\;{\int_{0}^{2^{k+1}}|\tilde{E}_*^{(q)}|^2dt\over
b_d(2^{k+1})}\right)\end{equation} is achieved. Suppose the maximum
will be effected by the first entry in (\ref{nuda}).

Then clearly \begin{equation}\label{nudas}{\int_0^\infty
|\tilde{E}^{(q)}|^2 dt \over b_d(q)}\;\lesssim\; {\int_0^\infty
|\tilde{E}^{(q)}|^2 dt \over b_d(2^{\bar{k}+1})},\end{equation} that
is $2^{\bar{k}+1}\lesssim q$. Then let $\bar{q} = 2^{\bar{k}+1}$,
and now look back at (\ref{sovsempochti}) with $\bar{q}$ instead of
$q$. Then by (\ref{ustal}) the terms in the right-hand side of
(\ref{sovsempochti}) corresponding to $k\lesssim \bar k$ would
change by at most a constant factor. On the other hand once more
$m(\bar{q})$ should be achieved for some $k\lesssim \bar{k}$. That
means that $m(\bar{q})\approx m(q)$, that is the supremum in
(\ref{sovsempochti}) with $\bar{q}$ substituting $q$ is of the same
order as the left-hand side, hence it is $O(1)$. By definition of
$\bar{q}$ this is thence the case for (\ref{sovsempochti}) per se as
well. This completes the proof of Lemma \ref{ebound} and Theorem 3.
$\Box$

\vspace{.2in} \noindent{\large {\bf Acknowledgements:}} Research has
been partially supported by the NSF Grant DMS02-45369, EPSRC Grant
GR/S13682/01 and the Bristol Institute for Advanced Study.

\medskip
\noindent {\large Alex Iosevich:} University of Missouri, Columbia,
Missouri,
65211 USA, {\em iosevich@math.missouri.edu}\\
{\large Mischa Rudnev:} University of Bristol, Bristol BS6 6AL UK,
{\em m.rudnev@bris.ac.uk}

\end{document}